\documentclass[12pt]{article}

\usepackage{amssymb}
\usepackage{amsthm,amsfonts,amsmath,amssymb}

\usepackage{amscd}

\newcommand{\eq}{\begin{equation}}
\newcommand{\en}{\end{equation}}

\newcommand{\diff}{{\tt d }}

\newcommand{\Nat}{\Bbb N}

\newcommand{\ed}{ \stackrel{d}{=}}

\newcommand{\stat}{{\tt stat }}

\newcommand{\rec}{{\tt rec}}

\def\endpf{\hfill $\Box$ \vskip0.5cm}

\newtheorem{theorem}{\large Theorem}
\newtheorem{proposition}[theorem] {\large Proposition}
\newtheorem{construction}[theorem] {\large Construction}

\newtheorem{lemma}[theorem]{\large Lemma}

\newcommand{\M}{{\rm M}}

\begin{document}

\title{Coherent random permutations with record statistics}
\author {Alexander Gnedin\thanks{Insitute of Mathematics, Utrecht University, The Netherlands, gnedin@math.uu.nl }}
\maketitle

\begin{abstract}\noindent
Random permutations with distribution conditionally uniform
given the set of record values can be 
generated 
in a unified way, coherently for all values of $n$. 
Our central example is  a two-parameter family of
random permutations that are conditionally uniform given the counts of upper and lower records.
This family interpolates between 
 two versions of Ewens' distribution. 
We discuss characterisations of the conditionally uniform permutations, 
their asymptotic properties, constructions and relations to random partitions.
\end{abstract}

\section{Introduction}

Random permutations $\pi_n\in {\frak S}_n$ with  distribution  conditionally   uniform
 given the value of some statistic 
${\tt stat}$
offer  a wide and most natural generalisation of the uniform distribution on the symmetric group ${\frak S}_n$. 
It is sometimes possible to define random permutations coherently for all 
values of $n$,  in a way  connecting the asymptotic properties of $\pi_n$'s  with a de Finetti-type representation
which generates $(\pi_n)$ from some limiting form of $\stat$ by means of a standard sampling procedure.

\par The most studied instance of coherent permutations is the one with $\stat$   defined as the 
nondecreasing sequence of cycle-sizes of $\pi_n$, see \cite{ABT, CSP}. In this case 
the sequence $({\stat}(\pi_n),~n=1,2,\ldots)$ 
is  Kingman's partition structure. The scaled cycle-sizes converge to a 
nonincreasing array
of frequencies $(p_k)$, from which $(\pi_n)$ can be recovered by a stochastic algorithm known as Kingman's
paintbox process. 
A distinguished example of coherent permutations with the cycle statistic 
is the parametric family of Ewens' distributions for $(\pi_n)$,  
associated  with the Poisson-Dirichlet law for the frequencies.
Ewens' distributions, as well as more general  two-parameter distributions  due to Pitman,
can be generated by  a simple urn scheme which does not exploit 
the  asymptotic frequencies \cite{CSP}.
In fact, for  Ewens' family   the minimal sufficient statistic is just the number of cycles of $\pi_n$
and, moreover, every distribution for $(\pi_n)$ with this property is a unique mixture
of Ewens' distributions, see \cite[Theorem 12 (i)]{Gibbs}.
See  \cite{ Gibbs, RPS, CSP} and references therein for many examples of partition structures.

\par Adopting for $\stat$ the sequence of cycle-sizes arranged by 
increase of minimal elements of the cycles leads to a wider type of structure 
introduced by Pitman \cite{PTRF} in the form of partially exchangeable partitions. 
The blocks of these ordered partitions correspond to the cycles  
of permutation, with 
the circular arrangement of 
elements within the cycles ignored.
In this case the limiting shape of $\stat$ 
is an arbitrary   random array $(p_k)$ of positive frequencies adding to at most unity.
For instance, for Ewens' permutations  this sequence
of frequencies has the GEM distribution, 
which is a size-biased arrangement 
of the Poisson-Dirichlet law.

\par Gnedin and Olshanski \cite{GO} studied coherent permutations with the set of descent positions of $\pi_n$ in the role of $\stat$.
They showed that coherent sequence of permutations $(\pi_n)$ corresponds to a spreadable  random  order on $\Nat$
(spreadability, also called contractability, means
invariance under all increasing injections  $\Nat\to\Nat$, see \cite{Kallenberg}).
The limiting shape of ${\stat}$ was identified with two disjoint open subsets 
of $[0,1]$. 
It was further shown in \cite{Euler} that if the law of
each $\pi_n$ is uniform conditionally given the number of descents, then $(\pi_n)$
is a unique mixture  of $a$-shuffles (introduced in \cite{shuffle}) and reversed $a$-shuffles.
In this sense the $a$-shuffles in the setting with descent statistic can be regarded as   analogues of Ewens' distributions
in the setting with cycle statistic. See \cite[Section 8.6]{GO} for results in the setting where $\stat$ is the peak set of permutation.

\par  Kerov and Tsilevich \cite{KT, Kerov} 
studied coherent permutations with $\stat$ defined to be the set of upper records of $\pi_n$.
This structure can be reduced 
to Pitman's \cite{PTRF}  partially exchangeable partitions  by the virtue of a fundamental
bijection ${\frak S}_n\to{\frak S}_n$ which translates 
the record statistics into the cycle statistics \cite[p. 17]{Stanley}.
In the interpretation in terms of records, 
the role of limiting shape of $\stat$ is played by partial sums $(p_1, p_1+p_2,\ldots)$,
which are also the upper record values of a random sequence $(X_n)$, 
such that the $\pi_n$'s can be generated  by ranking the variables $(X_n)$.
See \cite{CR, GnPEP} for more on partially exchangeable partitions and an application to multivariate records.

\par In this paper ${\stat}$ is the  two-sided set of  records of $\pi_n$, both upper and lower. 
We extend  known results \cite{GoldieRec, Kerov, KT, PTRF}
to include both types of records  in a symmetric way.
In particular,  
Ewens' family of random permutations  
will be extended  to a two-parameter  family of distributions $P^{(\theta,\zeta)}$.
Generalising the above mentioned one-sided result \cite[Theorem 12 (i)]{Gibbs} we  show that every coherent $(\pi_n)$ with $\pi_n$ conditionally uniform given 
the counts of upper and lower records is a mixture of the $P^{(\theta,\zeta)}$'s.
Permutations under $P^{(\theta,\zeta)}$ can be generated by ranking a sequence of real-valued random variables 
$(X_n)$, whose records follow a two-sided analogue of the GEM distribution.
This kind of representation  is also shown
for arbitrary coherent sequence of permutations $(\pi_n)$, with each $\pi_n$ uniformly distributed
given its set of record values.
Explicit formulas  are possible for a  multiparametric
class of distributions for $(\pi_n)$, 
which may be regarded as a  two-sided generalisation of a well-known Pitman's construction of exchangeable 
partitions \cite{PTRF}.

 \section{Counting the records}

Permutations $\pi_n\in{\frak S}_n$ of $[n]:=\{1,\ldots,n\}$ will be written 
in the one-row notation as 
$\pi_n=(\pi_{n1},\ldots,\pi_{nn})$. We call element
 $\pi_{nj}$ a  lower record  of $\pi_n$
if $\pi_{nj}=\min(\pi_{n1},\ldots,\pi_{nj})$, and we call
 $\pi_{nj}$ an upper record 
if $\pi_{nj}=\max(\pi_{n1},\ldots,\pi_{nj})$.
When $\pi_{nj}$ is a record we say that $\pi_{nj}$ is a record value and that
$j$ is  a record time (or a record position).
The first entry $\pi_{n1}$ will be called {\it center\,}. 
We regard the center as {\it improper}  lower and upper record, all other records being {\it proper}.
We denote  
$$\rec(\pi_n)=(r_{-\ell},\ldots,r_{-1}, r_0,r_1,\ldots, r_u)$$  
 the two-sided increasing sequence of record values, with distinguished center $r_0=\pi_{n1}$, proper lower records
$r_{-\ell},\ldots,r_{-1}$ and  proper upper records  $r_1,\ldots, r_u$.
In this notation  $\ell, u$ count the proper  records; 
for instance, 
$\rec(3,2,7,6,1,4,8,5)= (1,2,{\bf 3},7,8)$, where  the center is boldfaced and $\ell=u=2$.
 Clearly,  $r_{-\ell}=1, ~r_u=n$, and  the total number of records $\#\rec(\pi_n)=\ell+u+1$ 
satisfies $\min(2,n)\leq \ell+u+1\leq n$.
The record times of  proper lower and upper records will be labelled $t_1,\ldots,t_u$ and $t_{-1},\ldots,t_{-\ell}$,  respectively,
and we denote $t_0=1$ the record time associated with the improper record.

\par  Let $\left[\!{n\atop\ell+1,u+1}\!\right]$ be the number of permutations $\pi_n\in {\frak S}_n$ with
$\ell+1$ lower and $u+1$ upper records. This array of combinatorial numbers is symmetric in $\ell$ and $u$, and satisfies the
recursion
\eq\label{recS}
\left[\!{n\atop\ell+1,u+1}\!\right] =\left[\!{n-1\atop\ell ,u+1}\!\right]   +
\left[\!{n-1\atop\ell+1,u}\!\right]
+(n-2)
\left[\!{n-1\atop\ell+1,u+1}\!\right].
\en
Summing  over one of the parameters, say $u$, yields 
 a signless  Stirling number of the first kind
 $$\left[\!{n\atop \ell +1}\!\right]=\sum_{u=0}^{n-1}  \left[\!{n\atop\ell+1,u+1}\!\right],$$
equal to the number of permutations with $\ell+1$ lower records.
A more delicate connection to the Stirling numbers appears via   the identity
\eq\label{toSt}
\left[\!{n\atop \ell+1,u+1}\!\right]=\left[\!{n-1\atop \ell+u}\!\right]{\ell+u\choose\ell}
\en
found in \cite[p. 179]{DB}, where it was derived by manipulation with generating functions.

\par For our purposes it is important to introduce 
yet another encoding of permutation into the sequence of 
{\it initial ranks}
$$i_j:=\#\{k: ~k\leq j,\,\pi_{nk}\geq \pi_{nj}\},~~~~j\in [n].$$
The correspondence $\pi_n\mapsto(i_1,\ldots,i_n)$ is a well-known bijection 
 between ${\frak S}_n$ and $[1]\times[2]\times\cdots\times[n]$.
Note that
$\pi_{nj}$ is a lower record if $i_j=1$, and an upper record if $i_j=j$.

\par In terms of the initial ranks a bijective proof of (\ref{toSt}) is easily acquired.
To this end,
consider the mapping which sends $\pi_{n}\in{\frak S}_n$ to  $\pi_{n-1}'\in{\frak S}_{n-1}$ so that the initial ranks
 are transformed as $(i_1,\ldots,i_n)\mapsto (i_1',\ldots,i_{n-1}')$ where
$i_{j-1}'=i_j 1(i_j<j)$  for $2\leq j\leq n$. Each proper record of $\pi_n$  is mapped bijectively to
a lower record of $\pi_{n-1}'$, and the record counts satisfy $\ell(\pi_n)+u(\pi_n)=\ell(\pi_n')+1$.
It is easily seen that  $2^r$ permutations $\pi_n$ are mapped to the same $\pi_{n-1}'$ each time when $\ell(\pi_{n-1}')+1=r$, 
and of these  $\pi_n$ there are ${r\choose\ell}$ permutations  with $\ell$ proper lower records.
Because $\pi_{n-1}'$ with $r$ lower records can be chosen in $\left[\!{n-1\atop r}\!\right]$ ways, the identity (\ref{toSt}) follows.

\par When a probability distribution $P_n$ is specified on ${\frak S}_n$, we consider $\pi_n$ as a random variable.
In particular, $P_n^{(1,1)}(\pi_n)\equiv 1/n!$ is the  uniform distribution (indices will be explained in the next section). 
The characteristic feature of the uniform
distribution is that the initial ranks are independent, with
each $i_j$ being uniformly distributed on $[j]$. Giving a probabilistic interpretation to (\ref{toSt}) we have:

\begin{lemma}\label{exch} Under the  uniform distribution $P^{(1,1)}_n$ for $\pi_n$, 
conditionally given the record counts $(\ell,u)$ and given the positions 
occupied by $\ell+u$ proper records, all ${\ell+u \choose\ell}$ allocations  of $\ell$ lower records within these $\ell+u$ positions
are  equally likely.
\end{lemma}

\section{A two-parameter family of random permutations}\label{twopara}

We introduce next a two-parameter deformation of the uniform distribution,
for which $(\ell,u)$ is a sufficient statistic, meaning that  given the record counts the distribution of $\pi_n$ is uniform.
\begin{proposition}\label{pro1}
For arbitrary positive $\theta$ and $\zeta$ the formula
\eq\label{2par}
P_n^{(\theta,\zeta)}(\pi_n)={\theta^{\ell}\zeta^{u}\over (\theta+\zeta)_{n-1}}
\en
defines a  distribution on ${\frak S}_n$, which 
assigns the same probability to every permutation with $\ell+1$  lower  and $u+1$ upper records. 
\end{proposition}
\noindent

\par Proving this amounts to alternative definition of  
$P^{(\theta,\zeta)}_n$ as the probability distribution under which the initial ranks are independent
and  satisfy $i_1=1$ and for $j>1$
$$i_j= \begin{cases}1 {\rm~ ~w.p.~~} \theta/( \theta+\zeta+j-2),\\
j {\rm ~~w.p.~~} \zeta/(\theta+\zeta+j-2), \\
r {\rm ~~w.p.~~} 1/( \theta+\zeta+j-2)~~{\rm for ~~}r=2,\ldots,j-1
\end{cases}$$
(w.p.=with probability).
Multiplying these out
it is seen that (\ref{2par}) is the probability
of any sequence $(i_2,\ldots,i_n)$ where $i_j=1$ occurs $\ell$ times and $i_j=j$ occurs $u$ times.
Thus, $P_n^{(\theta,\zeta)}$ is obtained from $P_n^{(1,1)}$ by tilting 
the probabilities of extreme values of the initial ranks.

\par The fact that the probabilities in (\ref{2par}) add to  unity is also equivalent to the formula for the bivariate generating function 
\eq\label{2gf}
\sum_{\ell,u} \left[\!{n\atop \ell+1,u+1}\!\right]\theta^{\ell}\zeta^{u}=(\theta+\zeta)_{n-1},
\en
which dates back to  at least \cite{FS}. 
For $\zeta=1$ this specialises as the well-known formula 
$$\sum_{\ell=0}^{n-1} 
\left[\!{n\atop \ell+1}\!\right]\theta^{\ell-1}
=(\theta+1)_{n-1}$$
for the generating function
of Stirling numbers.

\par Recall that
{\it ranking} associates  with any sequence of distinct reals $x_1,\ldots,x_n$ a sequence of ranks 
$\pi_{nj}=\#\{i\leq n: ~x_i\leq x_j\}$, also called the {\it ranking permutation}. 
Ranking for the sequences with repetitions will be introduced in Section \ref{permrec}.

\noindent
{\bf Integer parameters.} 
For  integer $\theta,\zeta$ the distribution $P_n^{(\theta,\zeta)}$ can be obtained as a projection of the uniform distribution
$P_{n+d}^{(1,1)}$ on ${\frak S}_{n+d}$, where  $d=\theta+\zeta-2$. 
To ease notation, for the rest of this section the elements of permutation are written with one index.
\par Fix $(w_1,\ldots,w_{n+d})\in {\frak S}_{n+d}$. 
A sequence $(\pi_j',~j\in [n])$  (which is a permutation of $n$ integers $\{\theta,\ldots,n+\theta-1\}$) is uniquely defined by the condition that $\{\pi_1',\ldots,\pi_j'\}\subset \{w_1,\ldots,w_{d+j}\}$ is the
subset of integers whose ranks among  $\{w_1,\ldots,w_{d+j}\}$ 
are neither among top $\zeta-1$ ranks nor among bottom $\theta-1$ ranks.
Here is the inductive definition. 
Let $s_1,\ldots,s_{n+d}$ be the initial ranks of $w_1,\ldots,w_{n+d}$.
At step $1$ we define  $\pi_1'$ to be the element of rank $\theta$ among $w_1,\ldots,w_{d+1}$,
thus leaving $\zeta-1$ elements ranked above and $\theta-1$ ranked below $\pi_1'$.
At step $j$ the element $w_{d+j}$ is added, if $\theta\leq s_{d+j}\leq j+\theta-1$ 
then $\pi_j'=w_{d+j}$,
if $1\leq s_{d+j}\leq \theta-1$ then $\pi_j'$ is defined to be the element of rank $\theta$ among $w_1,\ldots,w_{d+j}$, 
and if $j+\theta\leq s_{d+j}\leq j+d$ then $\pi_j'$ is defined to be the element of rank $j+\theta-1$ among $w_1,\ldots,w_{d+j}$.
Understanding the second arrow in $(w_1,\ldots,w_{n+d})\mapsto (\pi_1',\ldots,\pi_n')\mapsto (\pi_1,\ldots,\pi_n)$ 
as the ranking operation, we have defined a projection $f_{n}^{(\theta,\zeta)}$
from ${\frak S}_{n+d}$ to ${\frak S}_n$.
\begin{proposition}\label{integ}
For positive integers $\theta,\,\zeta$ the mapping $f_{n}^{(\theta,\zeta)}$ 
sends the uniform distribution on ${\frak S}_{n+d}$ (where $d=\theta+\zeta-2$) to 
$P^{(\theta,\zeta)}_n$.
\end{proposition}
\proof
In the above, the initial ranks for  $(\pi_1,\ldots,\pi_n)$ and $(\pi_1',\ldots,\pi_n')$ are the same, and 
are given for $j=2,\ldots,n$ by
\begin{eqnarray*}
i_j=\begin{cases}~~~ 1,~~~~~~~~~~~~\,\,{\rm if~~} s_{j+d}\in [1,\theta],  \\
s_{j+d}-\theta+1,~~{\rm if~~} s_{j+d}\in [\theta+1,j+\theta-2],\\
~~~j, ~~~~~~~~~~~~~\,{\rm if~~}s_{j+d}\in [j+\theta-1,j+d].
\end{cases}
\end{eqnarray*}
For uniform permutation, $s_{j+d}$ is uniform on $[j+d]$ and these are independent, hence the $r_j$'s are independent  
with respective probabilities  $\theta/(n+d-2), \zeta/(n+d-2)$ for extreme ranks and  equal probabilities for  other values of $i_j$.
\endproof

For irrational $\theta$ or $\zeta$ the distribution $P_n^{(\theta,\zeta)}$ cannot be obtained
as a projection of a uniform distribution on some combinatorial object.

\section{Coherent permutations}\label{coh}

Our view of permutation is biased towards  the interpretation as order, rather than mapping. Orders can be obviously restricted
from larger sets to smaller. In this direction, we say that permutations $\pi_n$ and $\pi_m$, for $m\leq n$, are {\it coherent}
if they determine the same order on $[m]$.
A sequence $(\pi_n)$ of coherent permutations $\pi_n\in{\frak S}_n$ defines a strict order $\lhd$ on the infinite set $\Nat$:
$j\lhd i$ iff $\pi_{nj}<\pi_{ni}$ for all $n\geq\max(j,i)$.

\par Let $D_{nm}:{\frak S}_n\to{\frak S}_{m}$ ($n>m$) be the projection which cuts the last $n-m$ entries of $\pi_n$ and replaces 
the first $m$ entries $\pi_{n1},\ldots,\pi_{nm}$ by their ranking permutation.
The projection $D_{nm}$ is the same 
as  restricting orders from $[n]$ to $[m]$, hence the coherence means that $D_{nm}(\pi_n)=\pi_m$. 
The space of all orders on $\Nat$ has the structure of 
the projective limit ${\frak S}^\infty:=\lim\limits_{\longleftarrow} {\frak S}_n$. 
This space ${\frak S}^\infty$  should not be confused with 
the infinite symmetric group ${\frak S}_\infty$ (of bijections $\Nat\to\Nat$ that displace  only finitely many integers), 
which is
the inductive limit of finite symmetric groups ${\frak S}_\infty:=\lim\limits_{\longrightarrow} {\frak S}_n$.

\par In terms of the initial ranks, $D_{nm}:(i_1,\ldots,i_n)\mapsto(i_1,\ldots,i_m)$
is just the projection on the first $m$ coordinates. 
Every infinite sequence $(i_n)$ determines an order $\lhd$ on $\Nat$, in which 
$n$ is ranked $i_n$th within the set $[n]$. 
Therefore ${\frak S}^\infty$ can be identified with the infinite product
space $[1]\times[2]\times\ldots$ 
Endowed with the product topology, ${\frak S}^\infty$ is a metrisable totally disconnected Borel space. 
When a probability measure is defined on 
${\frak S}^\infty$
we view $(\pi_n)\in{\frak S}^\infty$
as a random coherent sequence of permutations, or a random order on $\Nat$.
By the measure extension theorem, distributions $P_n$ on ${\frak S}_n$, defined for every $n$, determine a unique distribution $P$ 
on ${\frak S}^\infty$ for a coherent sequence of permutations if and only if the $P_n$'s are compatible with projections.  
\par We denote $P^{(\theta,\zeta)}$  the measure on ${\frak S}^\infty$ under
which the initial ranks $i_1,i_2,\ldots$ are independent,
with distribution as in Section \ref{twopara}. The  distributions  $(P^{(\theta,\zeta)}_n, ~n=1,2,\ldots)$ introduced in Proposition \ref{2par} 
 are coherent 
 projections of $P^{(\theta,\zeta)}$.
\par For an order $\lhd$ on $\Nat$ we shall say that an upper (or lower) record occurs at time $n$ if $i_n=n$ (respectively, $i_n=1$).
Reversing the order  is an automorphism of ${\frak S}^\infty$, which is written as either
$\pi_{nj}\mapsto n-\pi_{nj}$ for $j\in [n],~ n\in \Nat$, 
or, via the initial ranks, as $i_n\mapsto n-i_n$ for $n\in \Nat$.
Clearly, reversing the order swaps the types of records, hence maps $P^{(\theta,\zeta)}$ to $P^{(\zeta,\theta)}$.

\vskip0.5cm
\noindent
{\bf Remark.} Except $D_n:=D_{n,n-1}$ there are two other useful projections $D_n',D_n'':{\frak S}_n\to {\frak S}_{n-1}$.
Projection $D_n'$ deletes $n$ in the one-row notation of $\pi_n$, and $D_{n}''$ deletes $n$ in the cycle notation of $\pi_n$.
The projective limit $\lim\limits_{\longleftarrow} ({\frak S}_n,D_n'')$ was introduced in the representation theory  
of ${\frak S}_\infty$ as the space of {\it virtual permutations} \cite{KOV}, and  
$D_n'$ was used in \cite{GO}.
The isomorphism of three kinds of projective limits is established by means of the commutative diagram
$$
\begin{CD}
 \pi_n @>>>\pi_n^{-1}  @>>> (\pi_n^{-1}){\widehat\, }\\
D_n @VVV   
D_n' @VVV
D_n'' @VVV
\\
\pi_{n-1} @>>>   \pi_{n-1}^{-1} @>>>  (\pi_{n-1}^{-1}){\widehat\,}
\end{CD}
$$
where $\pi_n^{-1}$ denotes the inverse permutation, and $\pi_n^{\widehat\,}$ denotes the fundamental bijection of ${\frak S}_n$
which translates the one-row notation of permutation into the cycle notation of another permutation by inserting  parentheses `)(' before
each proper lower record, e.g. $(3, 2,7,6, 1,4,8,5)^{\widehat\,}=(3)(2,7,6)(1,4,8,5)$ 
(Stanley \cite[p. 17]{Stanley} gives a slightly different version of the mapping).

\section{Specialisations}\label{spec}

 Some special values of the parameters $\theta,\zeta$ and some limits are worth mentioning.
We call distribution $P$ on ${\frak S}^\infty$ {\it degenerate} if $P_n(\pi_n)=0$ for some $n$ and some $\pi_n\in {\frak S}_n$.
All distributions $P^{(\theta,\zeta)}$ for $\theta,\zeta>0$ are nondegenerate.

\noindent
{\bf The uniform distribution.} The measure  $P^{(1,1)}$ may be called the uniform distribution on ${\frak S}^\infty$, 
since every $P^{(1,1)}_n$ is the uniform distribution on ${\frak S}_n$, with $P^{(1,1)}_n(\pi_n)\equiv 1/n!$
for every $\pi_n\in {\frak S}_n$.
The corresponding random  order $\lhd$ on $\Nat$ has the characteristic  property of  {\it exchangeability}, 
that is the law of $\lhd$ is invariant under the action of ${\frak S}_\infty$.
This order appears by ranking an iid sample $(X_n)$ from the uniform distribution on $[0,1]$
(or some other contunuous distribution on reals).
For fixed $n$ there are also other ways to link uniform $\pi_n$ to a sequence 
of $n$ random reals \cite{GoldieRec}. 

\noindent
{\bf Ewens' distributions $P^{(\theta,1)}$ and  $P^{(1,\zeta)}$.}
Ewens' distribution on ${\frak S}_n$ (also called  $\theta$-biased permutation, see \cite{ABT})
is the one which assigns probability
$\theta^{c-1}/(\theta+1)_{n-1}$, 
to every permutation with $c$ cycles.
The partition of $n$ comprised of cycle-sizes of $\pi_n$ follows then the Ewens sampling formula.

\par Suppose $\zeta=1$, so the probabilities (\ref{2par}) become $P^{(\theta,1)}_n(\pi_n)=\theta^\ell/(\theta+1)_{n-1}$
where $\ell+1$ is the number of lower records of $\pi_n$.
When $\pi_n$ follows $P^{(\theta,1)}_n$ then also $\pi_n^{-1}$, because 
$\ell(\pi_n)=\ell(\pi_n^{-1})$. 
To see this, draw permutation in two dimensions as a point scatter $\{(j,\pi_{nj}),{j\in [n]}\}$.
Observe that the records are those points which do not have other points south-west of them. 
Flip the picture about the diagonal to see that the property is preserved.
The inversion combined with the $~{\widehat\,}~$-mapping in Section \ref{coh} transforms 
the distribution in its  conventional `cycle form'. Therefore we still call
$P^{(\theta,1)}$ and  $P^{(1,\zeta)}$ Ewens' distributions 
(this viewpoint was suggested  in \cite{KT}).
\par 
By the same flipping argument, the sequence of lower record times $t_{-\ell},\ldots, t_{-1}, t_0$
coincides with the decreasing sequence of lower record values of the inverse permutation $\pi_n^{-1}$,
hence under $P^{(\theta,1)}$ we have further symmetry:
$(t_{-\ell},\ldots, t_{-1}, t_0)\ed (r_0, \ldots, r_{-1}, r_{-\ell})$.

\noindent
{\bf Distributions with equal parameters.} For $\theta=\zeta$ there is a symmetry between lower and upper records. 
For distributions $P^{(\theta,\theta)}_n(\pi_n)=\theta^{\ell+u}/(2\theta)_{n+1}$
the minimal sufficient statistic is the total number of records $\ell+u+1$. Given the value of this statistic, $\pi_n$ 
is uniformly distributed.

\noindent
{\bf Bernoulli pyramids $P^{(\infty p,\infty (1-p))}~$ $(0\leq p\leq 1)$.}  If $\theta, \zeta\to\infty$ but so that $\theta/(\theta+\zeta)\to p$, 
then under the limiting law the probability of $\pi_n$ 
is $p^\ell(1-p)^u$ provided $\ell+u=n-1$, and the probability is zero otherwise. 
Such $\pi_n$ has each $\pi_{nj}$ $(j>1)$ an upper record with probability $p$ and a lower record with probability $1-p$. 
Only  extreme initial ranks
 are possible, i.e $i_j\in \{1,j\}$. Such distributions were 
exploited in optimal stopping \cite{GK}.
One way to generate such permutation is to split $[n]$ by binomial variable at some integer
$v$, then let $\pi_1=v$ for the center and then riffle-shuffle $v+1,\ldots,n$ and $v-1,\ldots,1$ to obtain $\pi_{2n},\ldots,\pi_{nn}$.
In the cases $p=1$ (respectively,  $p=0$)  the distribution concentrates on the permutation $(n,\ldots,1)$ 
(respectively, $(1,\ldots,n)$).

\noindent
{\bf Degenerate  Ewens' permutations $P^{(\theta,0)},~P^{(0,\zeta)}$.} 
In the limiting case $\theta\to 0$ (but $\zeta>0$), the permutation has the form $\pi_n=(1,\pi_{n-1}')$, where
 $\pi_{n-1}'$ is a permutation of $\{2,\ldots,n\}$ which upon obvious identification has $P_{n-1}^{(1,\zeta)}$ distribution.
In the limiting case $\zeta\to 0$ (but $\theta>0$), the permutation has the form $\pi_n=(n,\pi_{n-1}')$, where
 $\pi_{n-1}'$ is a permutation of $[n-1]$ which has $P_{n-1}^{(\theta,1)}$ distribution.

\noindent
{\bf Permutations with only one proper record   $P^{(p0,(1-p)0)}~$ $(0\leq p\leq 1)$.}
When both $\theta,\zeta\to 0$ but so that $\theta/(\theta+\zeta)\to p$ for some $p\in [0,1]$, then 
the limit law of $\pi_n$ is that of $(\pi_{n1},\pi_{n2},\pi'_{n-2})$ where 
$(\pi_{n1},\pi_{n2})$ is either $(1,n)$ or $(n,1)$ with probability $p$ and $1-p$, respectively, while $\pi'_{n-2}$
is a uniform permutation of $\{2,\ldots,n-1\}$ independent of $(\pi_{n1},\pi_{n2})$.

\begin{proposition} The weak closure of the $(\theta,\zeta)$-family is comprised of nondegenerate ditributions
with $\theta>0$, $\zeta>0$, and of three degenerate types described above.
\end{proposition}
\proof This follows by considering $P_2^{(\theta,\zeta)}$ and $P_3^{(\theta,\zeta)}$. 
\endproof


\section{Characterisation  of mixtures}\label{mix}

We seek now for a two-parameter generalisation of
\cite[Theorem 12 (i)]{Gibbs}, that is we wish to characterise 
the distributions $P^{(\theta,\zeta)}$ as extreme 
points of a suitable  family of conditionally uniform distributions. 
The following lemma is helpful.
\begin{lemma}\label{JFin} 
Let $Q_1$ be the law of an independent $0$-$1$ sequence $B_1,B_2,\ldots$ with $B_n$ {\rm Bernoulli}$(1/n)$. 
Assume $Q$ is a distribution for $B_1,B_2,\ldots$ 
with the property that, for each $n$, the conditional law of $(B_1,\ldots,B_n)$ 
given $S_n:=B_1+\ldots+B_n$ and given $(B_m, ~m>n)$ 
under $Q$ is the same as under $Q_1$. Then $Q$  is a unique mixture of distributions
$Q_\eta$, $\eta\in [0,\infty]$, under which $B_1,B_2,\ldots$ are independent with $B_n$ {\rm Bernoulli}$(\eta/(n+\eta-1))$. 
\end{lemma}
\noindent
\proof
This can be concluded from either  \cite[p. 269]{PitmanFin} or \cite[Lemma 9]{Gibbs}.
The key issue is that 
 the convergence $S_n/\log n\to \eta$  holds under $Q_\eta$ almost surely.
\endproof
\par The first two assertions of the next proposition are equivalent to 
\cite[Theorem 12 (i)]{Gibbs} and included here for completeness of exposition.
\begin{proposition}\label{extp1} Suppose under $P$ the law of $\pi_n$  for every $n=1,2,\ldots$ 
is uniform conditionally given the value of a statistic $\stat$. 
Then the following assertions are true:
\begin{itemize}
\item[{\rm (i)}]  for $\stat=\ell$ distribution $P$ is a unique mixture of $P^{(\theta,1)}$ $(\theta\in [0,\infty[)$ and 
$P^{(1\infty,0\infty)}$,
\item[{\rm (ii)}]  for $\stat=u$ distribution $P$ is a unique mixture of $P^{(1,\zeta)}$ $(\zeta\in [0,\infty[\,)$ and
$P^{(0\infty,1\infty)}$,

\item[{\rm (iii)}]  for $\stat=\ell+u$ distribution $P$ is a unique mixture of 
$P^{(\theta,\theta)}$ $(\theta\in ]0,\infty[\,),$ 
$P^{({1\over 2}0,{1\over 2}0)}$ and $P^{({1\over 2}\infty,{1\over 2}\infty)}$,
\item[{\rm (iv)}] for  $\stat=(\ell,u)$ distribution $P$ is a unique mixture 
of  nondegenerate distributions $P^{(\theta,\zeta)}$ $(\theta,\zeta\in ]0,\infty[\,)$,
degenerate distributions
$P^{(\theta,0)}$ and  $P^{(0,\zeta)}$  $(\theta,\zeta\in ]0,\infty[\,)$,
and further degenerate distributions $P^{(1\cdot0,0\cdot 0)},$ $P^{(0\cdot0,1 \cdot 0)}$ and
$P^{(p \infty,(1-p)\infty)}$ $(p\in [0,1])$.
The degenerate distributions  do not enter provided $P_3>0$.
\end{itemize}
\end{proposition}
\proof
We need to show that the described distributions and only they are extreme.
Assuming $P$  extreme in the setting of  (iv),
the tail algebra $\cal F$ of the process $((\ell(\pi_n),u(\pi_n)),~n=1,2,\ldots)$ must be trivial.
Let $B_n=1(r_{n+1}\in\{1,n+1\})$ be the indicator of some record at position $n+1$. 
Under $P^{(1,1)}$ the law of $(B_1,B_2,\ldots)$ is $Q_2$, hence 
by  Lemma \ref{JFin}  and because $\lim S_n/\log n$ is $\cal F$-measurable
the law of $(B_n)$ under $P$ is the same as under  
$Q_{\eta}$ for some $\eta$. 
This says that records occur 
by a Bernoulli process, without specifying the types of  records.
If $\eta=0$ the situation is clear: there is only one proper record (for $n>1$) and  
$P^{(1\cdot0,0\cdot 0)},$ $P^{(0\cdot0,1 \cdot 0)}$ are the sole possibilities. 
Suppose $\eta\neq 0$.
 A key to recognise  how the records are classified in types is the exchangeability.
Let $I_k$ be the indicator of the event that the record at $(k+1)$st record time is a lower record.
Conditionally given $I_1+\ldots+I_k=\ell-1$    all values of the sequence $(I_1,\ldots,I_k)$ 
have the same probability $1/{k \choose \ell-1}$, because by Lemma \ref{exch} this is true under $P^{(1,1)}$
and by a simple stopping times argument. By de Finetti's theorem,
there exists a relative frequency of lower records, hence $\ell(\pi_n)/(\ell(\pi_n)+u(\pi_n))$ must converge almost surely.
But the limit of this ratio is $\cal F$-measurable hence constant, say $p$.
Appealing again to Lemma \ref{exch} we see that $(B_n)$ and $(I_k)$ are independent, hence 
the set of positions of lower records is the one obtained by independent 
thinning with probability $p$ of the occurences of $1$'s in $(B_n)$. Thus
$P=P^{(\theta,\zeta)}$ with $\theta=p\eta,~\zeta=(1-p)\eta$ (the instance $\eta=\infty$ is included).
Part (iii) is shown similarly, with the special feature that $p=1/2$. \endpf

\vskip0.5cm
\noindent
{\bf Remark.} To put the last result in the framework of \cite{Gibbs, Euler},
denote
$w_n(\ell,u)$ the probability for $\ell$ lower and $u$ upper proper records in $\pi_n$. 
By the rule of addition of probabilities we have 
\eq\label{dual}
w_n(\ell,u)=w_{n+1}(\ell+1,u)+w_{n+1}(\ell,u+1)+(n-1)w_{n+1}(\ell,u), ~~~~w_1(0,0)=1,
\en
 which is a recursion dual to (\ref{recS}). 
The set of nonnegative solutions to (\ref{dual}) is a convex compact set.
Proposition \ref{extp1}(iv) describes the set of extreme solutions to (\ref{dual}).
Interestingly, the set of extremes is not closed: each distribution
$P^{(p0,(1-p)0)}~$ with $0< p<1$ appears as a limit of some nondegenerate $P^{(\theta,\zeta)}$'s, but it is 
decomposable  as a mixture $P^{(p0,(1-p)0)}=
pP^{(1\cdot0,0\cdot0)}+(1-p)P^{(0\cdot0,1\cdot0)}~$.

\par A common approach to finding the extreme solutions of (\ref{dual}) is based on the analysis of asymptotic regimes for 
$\ell'= \ell'(n'), u'=u'(n')$ as $n'\to\infty$, which guarantee for all $n,\ell,u$  convergence of the ratios
\eq\label{rat}
\left[{n\atop \ell+1,u+1}\right]_{{n'\atop \ell'+1,u'+1}}{\bigg/}\left[{n'\atop \ell'+1,u'+1}\right],
\en
where the numerator is the number of permutations $\pi_{n'}$ of $[n']$ with record counts $(\ell',u')$ such that the restriction 
of $\pi_{n'}$ to $[n]$ has record counts $(\ell,u)$. 
Using a monotonicity argument, the things can be reversed to show that the convergence
$\ell'/\log n'$ and $u'/\log n'$ is necessary and sufficient for the convergence of the ratios (\ref{rat}) for all $n,u,\ell$.

\section{Some properties and asymptotics} 

As in the case of uniform distribution \cite{Nevzorov}, asymptotic 
properties (as $n\to\infty$) of record counts $\ell,u$ under $P^{(\theta,\zeta)}$ follow straightforwardly 
from the representation via independent initial ranks.
Thus, both mean and variance of $\ell$ are asymptotic to $\theta\log n$, 
and that of $u$ to $\zeta\log n$.  Jointly, $(\ell,u)$ converge in distribution to independent Gaussian variables.    
The point processes of scaled  record times $\{t_k/n:~k<0\}$, $\{t_k/n:~k>0\}$
converge to independent Poisson processes with intensities $\theta\diff t/t$, $\zeta\diff t/t$ (for $t\in [0,1]$), respectively.

\par The behaviour of each $\pi_{nj}$ under
$P^{(\theta,\zeta)}$ as $n$ varies is that of a process with exchangeable $0$-$1$ increments,
known as  P{\'o}lya's urn model. That is to say,
each sequence
 $(\pi_{nj}, ~n\geq j)$ is a nondecreasing inhomogeneous Markov chain on integers, which starts at some
random initial rank $\pi_{jj}=i_j$  at time $j$, 
and at time $n$ either jumps from some rank $\pi_{nj}=v$ to  $v+1$ with probability
$(v-1+\theta)/(n-2+\theta+\zeta)$, or otherwise remains at $v$.

\par The law of $\rec(\pi_n)$ can be expressed in terms of P{\'o}lya-Eggenberger distributions
$${\rm PE}_n^{(\theta,\zeta)}(r):={n-1\choose r-1}{(\theta)_{n-1}(\zeta)_{r-1}\over(\theta+\zeta)_{n-1}}\,~~~~~r\in [n].$$
The distribution of the center $r_{0}=\pi_{n1}$ is ${\rm PE}_n^{(\theta,\zeta)}$.
Conditionally given $r_{0}$,  
the lower and upper record sequences
are independent.
The sequence of lower records $r_{-1},\ldots,r_{-\ell}$ 
is a homogeneous decreasing Markov chain on integers which starts at $r_0$ 
and terminates at $1$, each time
descending from the generic $r$ to $r-d$ with probability
${\rm PE}_r^{(\theta,1)}(d)$.
In a similar way,
the sequence of upper records $r_{1},\ldots, r_{u}$ is a homogeneous increasing Markov chain on integers which starts at $r_0$
and terminates at $n$,
each time ascending from some $r$ to $r+d$ with probability
${\rm PE}_{n-r+1}^{(\zeta,1)}(d)$.

\par Asymptotics of the record values follow from well known properties of P{\'o}lya urns. 
Recall that 
 beta$(a,b)$ distribution with parameters $a>0,~b>0$  is the distribution on $[0,1]$ with
 density $x^{a-1}(1-x)^{b-1}/{\rm B}(a,b)$, where ${\rm B}(a,b)=\Gamma(a)\Gamma(b)/\Gamma(a+b)$.

\begin{proposition}\label{SLLN} As $n\to\infty$, under  $P^{(\theta,\zeta)}$ the scaled 
record values of $\pi_n$ converge,
$${r_k\over n}\to \rho_k\,~~{\rm a.s.}~~~(k\in{\mathbb Z}).$$
The distribution of $\rho_0$ is {\rm beta$(\theta,\zeta)$}. Given $\rho_0$ the sequences 
$(\rho_{k},k<0)$ and $(\rho_{k},k>0)$ are independent and 
representable as  
$$\rho_k=r_0  T_{k}T_{k+1}\cdots T_{-1} ~~~~(k<0),~~~~~
~\rho_k=1-(1-r_0)Z_1Z_2\cdots Z_{k}~~~~(k>0),$$
where
$T_k$'s are {\rm beta}$(\theta,1)$,  $Z_k$'s are {\rm beta}$(\zeta,1)$ and the
variables 
$\rho_0$, $T_k$  $(k<0)$ and $Z_k$  $(k>0)$  are all independent.
\end{proposition}

\par Let  $\cal S$ be the space of two-sided nondecreasing sequences $(x_k,~k\in {\mathbb Z})$, $x_k\in [0,1]$.
  We endow $\cal S$ with the product topology of  $\prod_{k=-\infty}^\infty [0,1]$.
Padding $\rec(\pi_n)$ by infinitely many $1$'s on the left and infinitely many $n$'s on the right, and scaling by $n$ makes
$n^{-1}\rec(\pi_n)$  a random element of ${\cal S}$ 
 $$n^{-1}\rec(\pi_n) = (\ldots,1/n,1/n, r_{-\ell}/n\ldots,r_{-1}/n,r_0/n,r_1/n,\ldots,r_u/n,1,1,\ldots).$$ 
Proposition \ref{SLLN} 
is a strong law of large numbers which
says that $n^{-1}\rec(\pi_n)$ converge in $\cal S$ almost surely to a limiting `shape' $(\rho_k)$.

\par Recall that GEM$(\theta$) distribution is the law of the sequence of gaps obtained by breaking $[0,1]$ at atoms of the Poisson point process
with intensity $\theta{\rm d}x/x$ $(x\in [0,1]$).
The decreasing sequence of atoms has the same distribution as
the sequence of `stick-breaking' products $D_1,D_1D_2,\ldots$, with
the  $D_j$'s being iid beta($\theta,1)$.
\par The two-sided sequence $(\rho_k,~k\in {\mathbb Z})$ is obtained in a similar
way, by splitting  $~[0,1]$ at $\rho_0$, and  further partitioning
the intervals $[0,\rho_0]$ and $[\rho_0,0]$  by two independent beta stick-breakings with  parameters $\theta$ and $\zeta$.
By analogy, the sequence of gaps $\rho_{k+1}-\rho_k,~ k\in{\mathbb Z}$, may be regarded as
a two-sided version of GEM distribution. 
\par Generalising the classical case of sampling from iid uniforms 
\cite[Proposition 4.11.2]{Resnick}, 
the distribution of the bivariate point process of upper (or lower) record values and durations follows 
from the spraying property of Poisson processes.
Thus,
given $\rho_0$ the  point processes $\{(\rho_{k}, t_{k+1}-t_k),~k\geq 0\}$ and  $\{(\rho_{k}, t_{k-1}-t_k),~k\leq 0\}$,
are independent Poisson,
with intensity measures $\zeta x^{j-1}\diff x$ on $[\rho_0,1]\times\Nat$ and $\theta (1-x)^{j-1}\diff x$ on $[0,\rho_0]\times\Nat$, respectively.
In particular, by the projection property of Poisson processes, given $\rho_0$
the conditional distribution of the number of pairs of neighbouring lower records $\#\{k\leq 0:  t_{k-1}-t_k=1\}$ 
 is Poisson($\theta\rho_0$)
(an equivalent result is shown in \cite[Corollary 3.1]{Holst} by computation of moments).

\section{Generating random permutations}

Under $P^{(\theta,\zeta)}$ not only the scaled record values converge (see Proposition \ref{SLLN}), 
but also scaled permutations
$(\pi_{nj}/n,\,j\in \Nat)$ converge almost surely to some random sequence $(X_j)\in [0,1]^\infty$.
In the case of uniform distribution $P^{(1,1)}$, the sequence $(X_j)$ is just iid uniform$[0,1]$,
and $(\pi_n)$ can be generated by ranking $(X_j)$.
Under any $P^{(\theta,\zeta)}$,  $(X_j)$ can be produced by a kind of shuffling of the sequences 
of record values
$(\rho_k,~ k\geq 0)$, $(\rho_k,~ k< 0)$ and
another independent sequence of uniform variables.
Here and henceforth, under shuffling of a few sequences 
 we understand a sequence which is comprised of terms of all these sequences
arranged in such a way that each of the sequences enters in its original order.

\begin{construction}\label{constr}{\rm
 Let $(W_n)$ be iid uniform$[0,1]$, independent of $(\rho_k)$.
We define a new sequence $(X_n)$ where some $W_n$'s are used, and some are replaced by 
$\rho_k$'s  which will appear as upper and lower record values.
Start with $X_1=\rho_1$. Suppose before step $n+1$ the values $\rho_{-\ell},\ldots, \rho_u$ have been included into $X_1,\ldots,X_n$;
then $\rho_u=\max(X_1,\ldots,X_n)$ and  $\rho_{-\ell}=\min(X_1,\ldots, X_n)$.
 At step $n+1$ we let $X_{n+1}=\rho_{u+1}$ if $\pi_{n+1}>\rho_u$, or $X_{n+1}=\rho_{-\ell-1}$ if $\pi_{n+1}< \rho_{-\ell}$, 
or  $X_{n+1}=\pi_{n+1}$ otherwise.
Define 
a coherent sequence of permutations $(\pi_n)$ by ranking  $(X_n)$.}
\end{construction}

\noindent
It is obvious that, given $(\rho_k)$, the  sequence $(X_n)$ resulting from the construction has the same law as iid uniform$[0,1]$ sequence conditioned 
on its two-sided sequence of  record values (see \cite{GoldieBunge} for the one-sided case of upper records).
This works for any $\theta,\zeta$ because
conditionally given $(\rho_k)$ 
the distribution of $(\pi_n)$ under any $P^{(\theta,\zeta)}$ is the same as under
the uniform distribution $P^{(1,1)}$.

\par For every fixed $n$
a similar procedure yields  uniform permutation
$\pi_n$ conditioned on $\rec(\pi_n)$.
Start with setting $\pi_{n1}=r_0$. At each step $j>1$ we will have $\pi_{n1},\ldots,\pi_{n,j-1}$ already determined, with
some maximum $\max(\pi_{n1},\ldots,\pi_{n,j-1})=r_{u'}$  and some minimum $\min(\pi_{n1},\ldots,\pi_{n,j-1})=r_{-\ell'}$.
At step $j\in \{2,\ldots,n\}$ a value $v$ is chosen uniformly at random from $[n]\setminus\{\pi_{n1},\ldots,\pi_{n,j-1}\}$. 
If $v<r_{-\ell'}$ let $\pi_{nj}=r_{-\ell'-1}$,  if $v>r_{u'}$ let $\pi_{nj}=r_{u'+1}$, and if 
$r_{-\ell'}<v<r_{u'}$ let $\pi_{nj}=v$. 
The sampled value $v$ is replaced each time $v$  breaks the last upper or lower  record.
In $n$ steps
the increasing sequences $(r_{-\ell},\ldots,r_{-1})$, $(r_1,\ldots,r_{u})$ are shuffled with other elements of $[n]$. 
It is intiutively clear and not hard to show that,
as $n$ becomes large, $n^{-1}\rec(\pi_n)=n^{-1}(\ldots,1,r_{-\ell},\ldots r_{-1},r_0,r_1,\ldots,r_u,n,\ldots)$ 
will converge in $\cal S$ to $(\rho_k)$.
This is just because 
 sampling from large finite sets will have nearly the same effect as independent uniform choices from $[0,1]$.

\par
Apparently, from the viewpoint of statistical theory of extremes the sequence 
$(X_n)$ is rather  exotic, as it is chosen just to simulate desired
 behaviour of records. 
This differs general $P^{(\theta,\zeta)}$
from  the uniform distribution $P^{(1,1)}$, 
when `injecting' some extrinsic $(\rho_k)$ is not at all necessary
since the uniform sample $(W_n)$ supplies automatically appropriate record values, so  $(X_n)\ed (W_n)$.
Still, in the case of integer parameters there is a simpler way to produce appropriate $(X_n)$ from a sequence of uniforms,
as parallels the construction of permutations in Proposition \ref{integ}.

\noindent
{\bf Integer values of the parameters.}
The idea is to 
assume some `prehistorical' sample of uniforms.
Suppose $\theta\geq 1, \zeta\geq 1$ are integers.
For $d=\theta+\zeta-2$ let $V_1,\ldots,V_d,W_1,W_2,\ldots$ be iid uniform$[0,1]$.
 At step $1$ choose $X_1$ as the value of rank $\theta$ among $V_1,\ldots,V_{d}, W_1$.
At each step $n$ we will have $\max(X_1,\ldots,X_n)$ equal to the $(n-\theta+1)$th order statistic 
in  $V_1,\ldots,V_{d}, W_1,\ldots,W_n$, and 
$\min(X_1,\ldots,X_n)$ equal to the $\theta$th order statistic 
in  $X_1,\ldots,X_{d}, W_1,\ldots,W_n$. 
If $W_{n+1}> \max(X_1,\ldots,X_n)$ we set 
$X_{n+1}$ equal 
 to the $(n+\theta-1)$th order statistic 
in  $V_1,\ldots,V_{d}, W_1,\ldots,W_n, W_{n+1}$,
if $W_{n+1}< \min(X_1,\ldots,X_n)$ we set 
$X_{n+1}$ equal 
 to the $\theta$th order statistic 
in  $V_1,\ldots,V_{d}, W_1,\ldots,W_n, W_{n+1}$,
and otherwise let $X_{n+1}=W_{n+1}$.
This works, since
there are always $\theta$ spacings below $\min(X_1,\ldots, X_n)$ and $\zeta$ spacings above  $\min(X_1,\ldots, X_n)$,
thus the resulting ranking is as in the proof of Proposition \ref{integ}.
\vskip0.5cm

\par The described process shows that, for integer $\theta\geq1,\zeta\geq 1$, Proposition \ref{SLLN} is a consequence
of  properties of the uniform order statistics.
 For all other values of  $\theta,\zeta$  
the result can be 
interpolated   from the integer case, because
the law of each $\pi_n$ is a rational function of the parameters 
of beta laws for $T_{k},Z_k$.

\section{Permutations with the $\rec$ statistic}\label{permrec}

For arbitrary choice of the distribution for
$(\rho_k)\in {\cal S}$, there is some random sequence $(X_n)$  
resulting from
Construction \ref{constr}, such that given $(\rho_k)$ the law of $(X_n)$ is the same as for
independent uniforms conditioned on the record values.
This suggests that arbitrary coherent $(\pi_n)$ with each $\pi_n$ uniform given $\rec(\pi_n)$ can be derived 
in this manner.
In general, however, $(\rho_k)$ may have repetitions, therefore we need to be careful 
with defining permutations by  ranking. 

\par We are only interested in the sequences of reals 
 $x_1,x_2,\ldots$ with the property that if $x_i=x_j$ for $i\neq j$ then
$x_j=\max(x_1,\ldots,x_j)$ or $x_j=\min(x_1,\ldots,x_j)$. This means that only record values can be repeated.
We shall define now an order $\lhd$ on $\Nat$.
Suppose first that $x_1\neq x_2$, then 
we set $i\lhd j$ if  either (a): $x_i<x_j$,   
or (b): $i<j$ and $x_i=x_j=\max(x_1,\ldots,x_j)$, or (c): $j<i$ and $x_i=x_j=\min(x_1,\ldots,x_i)$.
The rules (b) and (c) are inconsistent if the sequence starts with $m>1$ repetitions
$x_1=\ldots=x_m\neq x_{m+1}$, in this case all rules apply for $i,j>m$ and we just require that 
each $j\leq m$ be attributed  the initial rank either $1$ or $j$ by some extrinsic rule.
For $(X_n)$ derived by Construction \ref{constr} from arbitrary random $(\rho_k)\in {\cal S}$ and independent uniform $(W_j)$, we 
define coherent sequence of permutations $(\pi_n)$ by ranking $(X_n)$, with account of these rules for repetitions.

\par For instance, for constant sequence $\rho_k\equiv p$, we obtain $X_n\equiv p$, and $i_n=1$ or $i_n=n$ according as $W_j<p$ or $W_j>p$,
so this $(\pi_n)$ is the  Bernoulli pyramid $P^{(p\cdot\infty,(1-p)\cdot\infty)}$.
Another example: permutations with single proper record,  $P^{(0\cdot p,0\cdot(1-p))}$, correspond to the case when $(\rho_0,\rho_1)=(0,1)$ w.p. $p$
and $(\rho_{-1},\rho_0)=(0,1)$ w.p. $1-p$.
Conditioning on $(\rho_k)$ and on $\rec(\pi_n)$ we have each $\pi_n$ uniformly distributed, whichever the values of
 $(\rho_k)$.

\par The main result says that 
this construction is indeed the most general.

\begin{proposition}\label{repre} Let $P$ be a distribution for  a coherent sequence of permutations 
$(\pi_n)$
with the property that, for every $n$, conditionally given $\rec(\pi_n)$, $P_n$ is a uniform distribution. 
Then $\rec(\pi_n)=(r_k)$ satisfies 

\eq\label{rlim}
{r_k\over n}\to \rho_k ~~~{\rm a.s.} ~~(k\in {\mathbb Z})
\en
for some random sequence $(\rho_k)$ with values in $\cal S$. 
Conditionally 
given  $(\rho_k)$, the law of $(\pi_n)$ is the same as for the coherent 
sequence of permutations generated by ranking the variables $(X_n)$ determined in 
{\rm Construction \ref{constr}}.
\end{proposition}

\par This is a de Finetti-type representation of $(\pi_n)$: given $(\rho_k)$,
the limit shape of the sufficient statistic $n^{-1}\rec(\pi_n)$, coherent permutations are generated by sampling uniforms
and shuffling them with $(\rho_k)$. 

\par One proof appeals to de Finetti's theorem for 0-1-sequences,  
and exploits the fact that given $n$ is the $k$th record time $t_k$ (so $i_n$ equals $1$ or $n$), the indicator variables $1(m\lhd n)$  for $m>n$
are exchangeable, where $\lhd$ is the order on $\Nat$ associated with $(\pi_n)$.
The exchangeability implies the existence of limits $(\ref{rlim})$.

\par Another proof is by reduction to Pitman's characterisation of partially exchangeable partitions \cite[Theorem 6]{PTRF}.
To this end, we need to associate with $(\pi_n)$ (thought of as order $\lhd$ on $\Nat$) an ordered partition $\Pi$ of $\Nat$ in disjoint nonempty blocks
$(A_k,\, k\in \mathbb Z)$. Let $A_0:=\{1\}$ be  
singleton block. For $k>0$ we assign to $A_k$ the $k$th proper upper  record time and all integers
$n$ $\lhd$-ordered between the $(k-1)$st and the $k$th proper upper record times.   Similarly, for $k<0$ 
we assign to  $A_k$  the $-k$th proper lower record time and all integers
$n$ $\lhd$-ordered between the $-k$th and the $(-k+1)$st proper lower record times. 
Thus the minimal elements
of blocks  are the record times $(t_k, ~k\in {\mathbb Z})$. We order the set of blocks $\{A_k,~k\in {\mathbb Z}\}$ 
by increase of the record values.
The sequences  $(t_k,~ k\geq0 )$ and $(t_k,~k\leq 0)$ start with common element $t_0=1$, are increasing and shuffled, that is 
interlaced in some random succession. 
Conditioning on the succession of record times $(t_k, ~k\in{\mathbb Z})$ (which could start like e.g. 
 $t_0,t_1,t_2,t_{-1}, t_3, t_{-2},\ldots$) we obtain a partially exchangeable partition, hence \cite[Theorem 6]{PTRF} 
can be applied, from which Proposition \ref{repre} follows by unconditioning.
\par The  differences $(\rho_{k+1}-\rho_k, ~k\in {\mathbb Z})$ are the frequencies $(p_k)$ 
of blocks of the ordered partition $\Pi$. In the event $\sup \rho_k<1$ or  $\inf \rho_k>0$
we have $\sum_{k\in {\mathbb Z}} p_k<1$ and  $\ell+u\sim (1-\sum_{k\in {\mathbb Z}} p_k)n$, i.e. the number of records grows linearly with $n$.

\section{The boundary of a composition poset}

The classification of coherent permutations with $\rec$ statistic fits in the Kerov-Vershik framework
of potential theory on graded graphs \cite{VK}.  We sketch this aspect of Proposition \ref{repre}.

\par Recall that $\rec(\pi_n)$ assumes values in the set of increasing sequences 
$r_{-\ell}<\ldots<r_0<\ldots<r_u$
 with the first term $1$, last term $n$ and
a distinguished center $r_0$. 
By a suitable differencing, $(r_k)$ can be bijectively encoded into
 a {\it centered composition} of integer $n$, which we define as  
 a sequence of positive integer parts 
$\lambda=(\lambda_{-\ell},\ldots,\lambda_{-1},\lambda_0,\lambda_1,\ldots,\lambda_u)$ 
with
distinguished {\it center}  $\lambda_0=1$ and $\sum_{k=-\ell}^u\lambda_k=n$.
The connection is established by the formulas
\begin{eqnarray}\label{connect}
r_k=1+\lambda_{-\ell}+\ldots+\lambda_{k-1}~~~(k\leq 0),~~~r_k=\lambda_{-\ell}+\ldots+\lambda_j~~~~(k>0),\\
\label{connect1}
\lambda_k=r_k-r_{k-1} ~~(k>0),~~~
\lambda_0=1,~~~
\lambda_k=r_{k+1}-r_k, ~~(k<0).
\end{eqnarray}

The centered composition $\lambda$ corresponding to $\rec(\pi_n)$
is the sequence of block-sizes of the 
ordered partition $\Pi_n=\Pi|_{[n]}$ from the previous section.

\par The number of centered
compositions of $n$ is $2^{n-3}(n+2)$. 
The centered compositions 
comprise a graded poset ${\cal R}$, in which immediate followers of $\lambda$ (centered composition of some $n$) are centered compositions $\mu$ obtained by either
incrementing one of noncentral parts by $1$ or by appending $1$ to the left or to the right.
For instance, $(3,1,{\bf1},3,2)$ is followed by $(1,3,1,{\bf1},3,2)$ ,
$(4,1,{\bf1},3,2)$, $(3,2,{\bf1},3,2)$, $(3,1,{\bf1},4,2)$, $(3,1,{\bf1},3,3)$ and $(3,1,{\bf1},3,2,1)$. 
The representation of increasing sequences by compositions is convenient because  passing to a follower requires incrementing only one part.

\par The {\it boundary problem} for $\cal C$ asks one to find all 
{\it extreme}
nonnegative solutions to the recursion 
$\phi(\lambda)=\sum_{\mu} \phi(\mu)$ with initial condition $\phi(1)=1$, where the summation is over $\mu$ which are immediate followers of $\lambda$.
The set of all nonnegative solutions is a compact convex set with the property 
that every its point has a unique representation as convex mixture of the extremes
(Choquet simplex).

\par By some well known general theory each extreme solution appears as a pointwise limit 
$\phi(\lambda)=\lim_{m\to\infty} d(\lambda,\mu)/d(\mu)$ for some sequence of centered compositions $\mu\in {\cal C}$ of growing 
degree $m\to\infty$. 
Here, $d(\lambda)$ 
is the number of permutations with $\rec(\pi_n)=(r_k)$ and $(r_k)$ corresponding to $\lambda$ via (\ref{connect}),
and $d(\lambda,\mu)$ is the number of permutations $\pi_m$ which correspond to $\mu$ 
and are coherent with some fixed permutation $\pi_n$ having this $\rec(\pi_n)=(r_k)$.
In other words, $d(\lambda,\mu)$ is the number of saturated chains in $\cal C$ which interpolate between centered compositions
$\lambda$ and $\mu$.

\par Computing the number   of permutations with fixed $\rec(\pi_n)$ yields
\eq\label{d1}
d(\lambda)=
{(n-1)!\over 
\Lambda_{-\ell}\cdots\Lambda_{-2}\Lambda_{-1}\Lambda_1\Lambda_2\cdots\Lambda_u}
\en
where
 $\Lambda_k=\lambda_k+\lambda_{k+1}\ldots+\lambda_u$ for $k>0$ and $\Lambda_k=\lambda_k+\lambda_{k-1}\ldots+\lambda_{-\ell}$ for $k<0$
are the right and the left tail-sums of $\lambda$.
 
If $\mu$ succeeds  $\lambda$ in $\cal C$ then $\mu$ is of the form
$\mu=(\mu_{-\ell-b},\ldots,\mu_{u+a})$ (for some $a,b\geq 0$) and $\mu_k\geq\lambda_k$ for $-\ell\leq k\leq u$.
For the number of $\pi_m$ coherent with $\pi_n$ we have 
\begin{eqnarray}\label{d2}
d(\lambda,\mu)  =
{(m-n)!\prod_{k=-\ell}^u{\mu_k-1\choose\lambda_k-1}\over
\M_{-\ell-b}\cdots\M_{-\ell-1}\M_{u+1}\cdots \M_{u+a}},
\end{eqnarray}
where $\M_k$ are tail-sums of $\mu$.
\footnote{Kerov \cite[Equations 1.4.4, 1.4.4]{Kerov} derived similar one-sided 
formulas from Stanley's dimension formula for coideals in trees. 
The factors  ${1\over n_j}{n_j\choose m_j}$ in \cite[Equation (1.4.4)]{Kerov} 
should be corrected as ${n_j-1\choose m_j-1}$.
The method of \cite{Kerov}  also 
applies here for a suitable tree and can be used to check (\ref{d1}), (\ref{d2}).}
From (\ref{d1}) and (\ref{d2}) 

\begin{eqnarray*}
d(\lambda,\mu)/d(\mu) =~
{\M_{-\ell}\cdots \M_{-1}\M_{1}\cdots\M_u\over (m-1)_{(n-1)\downarrow}}
\prod_{k=-\ell}^u {(\mu_k-1)_{(\lambda_k-1)\downarrow}\over (\lambda_k-1)!}\,,
\end{eqnarray*}
where $(x)_{k\downarrow}=x(x-1)\ldots(x-k+1)$ with $(x)_{0\downarrow}\equiv 1$.
Analysis of these explicit formulas shows that the ratios $d(\lambda,\mu)/d(\mu)$ converge as $m\to\infty$ for every $\lambda\in {\cal C}$
if and only if there exist limits $\mu_k/m\to p_k$ for each $k\in {\mathbb Z}$.  
In such limiting regime for $\mu$ with ${\bf p}:=(p_k)$, the resulting solution is 
\eq\label{poly}
\phi_{\bf p}(\lambda)= 
\prod_{k=-\ell}^{-1} \rho_{k+1} p_k^{\lambda_k-1}
\prod_{k=1}^u(1-\rho_{k-1})p_k^{\lambda_k-1}
\en
where $\rho_k=\sum_{i=-\infty}^k p_i$. 
Note that for any $\pi_n\in {\frak S}_n$ with $\rec(\pi_n)=(r_k)$ (corresponding to $\lambda$) we have
$\phi_{\bf p}(\lambda)=P_n(\pi_n)$ where $P=(P_n)$ is the distribution derived from $(\rho_k)$ by Construction \ref{constr}.
From the law of large numbers for this $P$ now follows that each $\phi_{\bf p}$ is an extreme solution.
This again implies 
Proposition \ref{repre}.

\par Finally, we mention one algebraic aspect. For each fixed $\lambda\in{\cal C}$ consider $\phi_{\bullet}(\lambda)$
as a formal polynomial (\ref{poly}) in infinitely many variables $(p_k,~k\in{\mathbb Z})$.
For various $\lambda\in{\cal C}$ these polynomials form a basis of an algebra $\cal A$, which has the 
property that the structural constants of multiplication in this basis are all nonnegative.
Moreover, $(\sum_{k\in {\mathbb Z}}p_k)\phi_{\bf p}(\lambda)=\sum_\mu \phi_{\bf p}(\lambda)$, where the sum is over immediate
followers $\mu$ of $\lambda$. In terms of \cite{VK} this means that the graded poset $\cal C$ is multiplicative.
By the Kerov-Vershik ring theorem (see \cite[Section 8.7]{GO} for detailed proof) extreme solutions have the form
$\phi(\lambda)=\chi(\phi_{\bullet}(\lambda))$
where $\chi$ is a 
homomorphism
$\chi:{\cal A}\to{\mathbb R}$ of algebras,
which satisfies $\chi \left(\sum_{k\in{\mathbb Z}}p_k \right)=1$ and also satisfies the positivity condition
$\chi(\phi_{\bullet}(\lambda))\geq 0$ for 
$\lambda\in {\cal C}$. Proposition \ref{repre} parametrises  all such $\chi$ by sequences
$(\rho_k)\in{\cal S}$, so that on the basis 
$(\phi_{\bullet}(\lambda), \lambda\in {\cal C})$ the homorphism 
is the specialisation $\phi_{\bullet}(\lambda)\mapsto \phi_{\bf p}(\lambda)$
with ${\bf p}=(p_k)$, where
$p_k=\rho_{k}-\rho_{k-1}$ for $k>0$ and $p_k=\rho_{k}-\rho_{k+1}$ for $k<0$.

\section{Further examples}

P{\'o}lya's urns allow to construct a large family of distributions for $(\pi_n)$ that 
are analogous to Pitman's two-parameter partition structures.
The idea is to extend the construction of the `Chinese restaurant process' \cite[Section 3.2]{CSP}
by tilting probabilities of extreme ranks together with intermediate ranks. 
\par Let $\alpha_k$ $(k\in{\mathbb Z})$, $\theta,\zeta\in{\mathbb R}$ be parameters.
Consider distribution $P$ such that given 
$i_1,\ldots,i_n$ the next initial rank satisfies

\eq\label{largef}
i_{n+1}=\begin{cases}
1~~~~~~~{\rm w.p.~~~~} {\theta+\alpha_{-1}+\ldots+\alpha_{-\ell}\over \theta+\zeta+n-1}\\
r~~~~~~~{\rm w.p.~~~~~~} {1-\alpha_{k-1}\over \theta+\zeta+n-1}~~~~~~~{\rm for~~} r_{k-1}<r\leq r_k,~~k<0\\
r~~~~~~~{\rm w.p.~~~~~~} {1-\alpha_k\over \theta+\zeta+n-1}~~~~~~~{\rm for~~} r_{k-1}<r\leq r_k,~~k>0\\
n+1~~{\rm w.p.~~~~} {\zeta+\alpha_1+\ldots+\alpha_{u}\over \theta+\zeta+n-1}
\end{cases}
\en
The {\it principal domain} of parameters is defined by the conditions of strict positivity
$$1-\alpha_k>0~~(k\neq 0), ~~~\theta+\alpha_{-1}+\ldots+\alpha_{-\ell}>0~~(\ell\in\Nat),~~~
\zeta+\alpha_1+\ldots \alpha_u>0~~(u\in\Nat).$$
Parameter $\alpha_0$ can be selected arbitrarily.
Under such $P$ the probability of every permutation $\pi_n\in {\frak S}_n$ with $\rec(\pi_n)=(r_k)$ is 
\begin{eqnarray*}
\phi(\lambda_{-\ell},\ldots,1,\ldots,\lambda_u)=~~~~~~~~~~~~~~~~~~~~~~~~~~~~~~~~~~~~~~~~~~~~~~~~~~~~~~~~~~~~~~~~~~~~~~~~~~~~~~~~~~~~~~~~~~~~~~~~~\\
{(\theta+\alpha_{-1})(\theta+\alpha_{-1}+\alpha_{-2})\ldots (\theta+\alpha_{-1}+\cdots+\alpha_{-\ell})
(\zeta+\alpha_{1})(\zeta+\alpha_{1}+\alpha_{2})\cdots (\zeta+\alpha_{1}+\cdots+\alpha_{u})
\over (\theta+\zeta)_{n-1}}\times \\
\prod_{k=-\ell}^u (1-\alpha_k)_{\lambda_k-1}\,,~~~~~~~~~~~~~~~~~~~~~~~~~~~~~~~~~~~~~~~~~~~~~~~~~~~~~~~~~~~~~~~~~~~~~~~~~~~~~~~~~~~~~~~~~~~~~~~~~~~
\end{eqnarray*}
where the centered composition $\lambda=(\lambda_{-\ell},\ldots,1,\ldots,\lambda_u)$ encodes $(r_k)$ via (\ref{connect}),(\ref{connect1}).
For the parameters in the principal domain the coherent permutations $(\pi_n)$ are nondegenerate.

The instance $\alpha_k\equiv 0$ corresponds to the $P^{(\theta,\zeta)}$-family.
Generalising Proposition \ref{SLLN} and specialising Proposition \ref{repre} we have the following representation.

\begin{proposition}\label{SLLN1} Suppose $P$ is defined by the conditional distributions {\rm (\ref{largef})}, with parameters in the principal
domain. Then under  $P$
   the scaled 
record values of $\pi_n$ converge, as $n\to\infty$,
$${r_k\over n}\to \rho_k\,~~{\rm a.s.}~~~(k\in{\mathbb Z}).$$
The distribution of $\rho_0$ is {\rm beta$(\theta,\zeta)$}. Given $\rho_0$ the sequences 
$(\rho_{k},k<0)$ and $(\rho_{k},k>0)$ are independent and 
representable as  
$$\rho_k=r_0  T_{k}T_{k+1}\cdots T_{-1} ~~~~(k<0),~~~~~
~\rho_k=1-(1-r_0)Z_1Z_2\cdots Z_{k}~~~~(k>0),$$
where
$T_k$'s are {\rm beta}$(\theta+\alpha_{-k}+\alpha_{-k+1}+\ldots+\alpha_{-1},1-\alpha_k)$,  $Z_k$'s are {\rm beta}$(\zeta+\alpha_1+\ldots+\alpha_k,1-\alpha_k)$ and the
variables 
$\rho_0$, $T_k$  $(k<0)$ and $Z_k$  $(k>0)$  are all independent.
\end{proposition} 
\noindent

\par Asymptotic properties of $\pi_n$ depend essentially on the parameters.
For instance, if $\alpha_k=a\in \,]0,1[$ for all $k>0$ and $\alpha_k=b\in\, ]0,1[$ for all $k<0$, then the order of growth of the number of upper 
records is $n^a$, and of the number of lower records is $n^b$, very much in line with asymptotics of 
Pitman's partitions \cite[Section 3.3]{CSP}. Extensions for other values of parameters, 
including those outside the principal domain, seem to be unexplored even in the one-sided case of  
upper records as sufficient statistic
(or partially exchangeable partitions).

\end{document}